\documentclass{amsart}

\newcommand\myshorttitle{Hitting time and fluid approximation: application
  to the coupon collector problem. }

\usepackage[utf8]{inputenc}
\usepackage[american]{babel}
\usepackage{url}
\usepackage{amsmath,amsfonts}
\usepackage{subfigure}
\usepackage{tikz}
\graphicspath{{fig/}{}}
\usepackage{amsthm}
\usepackage[colorlinks=true]{hyperref}

\newcommand\expect[1]{\mathbb{E}\left[#1\right]}
\newcommand\proba[1]{\mathbb{P}\left(#1\right)}
\newcommand\Sp{\mathcal{S}} 
\newcommand\F{\mathcal{F}} 
\newcommand\N{\mathbb{N}}
\newcommand\R{\mathbb{R}}
\newcommand\ceil[1]{\lceil#1\rceil}

\newcommand\norm[1]{\left\|#1\right\|}
\newcommand\lp{\left(}
\newcommand\rp{\right)}
\newcommand\p[1]{\left(#1\right)}
\newcommand{\bydef}{\stackrel{\rm{def}}{=}}

\newtheorem{theorem}{Theorem}
\newtheorem{lemma}{Lemma}

\begin{document}

\title[\myshorttitle]{Computing hitting times via fluid approximation:
  application to the coupon collector problem}

\author{Nicolas Gast} %
\thanks{Contact: \url{nicolas.gast@epfl.ch}.}%
\address{Nicolas Gast, EPFL, IC-LCA2, BC 203 Bâtiment BC, Station 14, 1015
  Lausanne-EPFL, Switzerland}%
\email{\url{nicolas.gast@epfl.fr}}

\begin{abstract}
  In this paper, we show how to use stochastic approximation to compute
  hitting time of a stochastic process, based on the study of the time for
  a fluid approximation of this process to be at distance $1/N$ of its
  fixed point.

  This approach is developed to study a generalized version of the coupon
  collector problem. The system is composed by $N$ independent identical
  Markov chains. At each time step, one Markov chain is picked at random
  and performs one transition. We show that the time at which all chains
  have hit the same state is bounded by $c_1N\log N + c_2 N\log\log N +
  O(N)$ where $c_1$ and $c_2$ are two constants depending on eigenvalues of
  the Markov chain.
  
\end{abstract}

\maketitle


\section{Introduction}

The coupon collector is a classical problem in probability theory. There
are $N$ types of coupons. Coupon are collected at random with
replacement. The goal is to compute the number of coupons to be collected
to have at least one coupon of each kind. This simple problem has a simple
answer and a simple proof: on average, one has to buy
$N/N+(N-1)/N+\dots+1/N\approx N\log N$ coupons to complete a collection: if
you already have $k$ different types of coupons, it takes in average
$(N-k)/N$ to get a coupon of a new type.

Because of its simplicity, this problem has many applications, especially
in computer science where it often serves as a basic tool for computing the
completion time of randomized algorithm
\cite{mitzenmacher2005probability,kenthapadi2005decentralized}.  Many
variants of it have been studied during the years. For example, the time
needed to complete a $T$ collections of the same $N$ coupons is shown to be
$N(\log N+(T-1)\log\log N+O(1))$ in \cite{double-dixie,myers2003some}. The
time to complete the first collection is $N\log N$. However, the time to
complete each next collection is only $N \log \log N$.

However, even a slight modification such has obtaining $T$ different
collections instead of one leads to much more complicated proofs. The
approach taken in this paper aims at being more general but also at giving
a new insight on the relation between hitting time and stochastic
approximation.

\subsection*{Contributions}

We develop an approach based on stochastic approximation to compute the
hitting time of a stochastic process that has an absorbing state.  The
system is composed of $N$ identical Markov chains that have an absorbing
state $0$.  At each time step, one chain is picked at random and performs
one transition. Our goal is to compute the number of steps until which all
Markov chains are in their absorbing state. The coupon collector problem is
a particular case of this problem by considering $N$ deterministic Markov
chains that have state $1$ (no coupon of that type has been collected) or
$0$ (at least one coupon of type $i$ has been collected).

Using a classical stochastic approximation approach, like
\cite{benaimLeboudec}, one can show that if $N$ is large, the proportion of
Markov chains that are in a given states can be approximated by a linear
ordinary differential equation (ODE) $\dot{m}=mQ$. This ODE has a unique
fixed point to which all trajectories converge exponentially fast,
corresponding to a state where all chains are in state $0$. However, this
approximation is not enough accurate to bound the hitting time of the
stochastic process. The time for this ODE to reach its equilibrium is
infinite whereas the expected time for the stochastic system to hit this
equilibrium is finite. 

In this paper, we establish a relation between the expected time
$\expect{T_N}$ for all the chains to be completed and the time $t_N$ for
the ODE to be at distance $1/N$ of its equilibrium point.  The main results
of this paper are Theorem~\ref{th:ODE_absorb} and
Theorem~\ref{th:bound_log}. We first show that $T_N$ is bounded by $N\cdot
t_N$. Using this result, we derive the existence of constants $c_1,c_2$
that depend on the spectral properties of the original Markov chain such
that
\begin{align*}
  \expect{T_N} &\le Nt_N + O(N) \\
  &\le c_1N\log N + c_2 N \log \log N + O(N).
\end{align*}
Applied to the time to complete $T$ collections, the allows to derive
directly the results of \cite{double-dixie}.

We also study two particular cases for which we have simple closed-form
bounds for $T_N$. In the more general case, if we only know that the
expected hitting time of one Markov chain is bounded by $T$ starting from
its initial state, then $T_N$ is only bounded by $N^2T$. However, if the
expected hitting for a Markov chain is bounded by $T$ independently of its
initial state, then we show that $\expect{T_N}\le NT \log N + O(N)$.  We
provide examples that show that these bounds are tight up to a linear term.
Finally, we also show that this method can be applied to study the
completion time of distributed algorithm.

We believe that the interest of this method is twofold. On the one hand, it
gives a new insight on the coupon collector by providing a new proof of a
more general results. We also think that this results could be adapted to
more general stochastic approximation algorithms and that this could be
helpful to understand the relationship between the extinction time of
stochastic models and the time for a fluid approximation of it to get close
to extinction.

\subsection*{Related work}

Stochastic approximation algorithm have been introduced in
\cite{robbins1951stochastic} for solving root finding problems. Application
of these methods are scattered on many fields, like economics
\cite{benaim-weibull} or computer science \cite{benaimLeboudec}. In all of
these works, a first step is to show that the stochastic system can be
approximated on any finite time interval by a fluid approximation --
\emph{e.g.} described by a differential equation. Then, this approximation
is used to derive asymptotic properties such as characterizing the limiting
dynamics \cite{benaim1999dynamics}, computing approximation of the
steady-state distribution \cite{benaim-weibull, benaimLeboudec} or proving
stability properties \cite{dai1995stability,fort2008ode}.

However, there are few results on the relation between the time for a
stochastic process to escape a region and the behavior of a fluid
approximation of it. In \cite{darling}, the authors shows that
if the time for the differential equation to escape a region is finite,
then the time for the stochastic system to escape this region converges to
the same value (Theorem~4.3 of \cite{darling}). When the differential
equation stays far from the absorbing boundary, the time for the stochastic
system to reach its absorbing state can be bounded by large-deviation
results, as in \cite{klebaner2001asymptotic}. Our case of interest in this
paper is that the deterministic system converges to a fixed point but does
not reach this point is finite time while the stochastic process does hit
this point in finite time.

\subsection*{Outline of the paper}

The rest of the paper is organized as follows. In
Section~\ref{sec:notations}, we give a formal definition of the
problem. Section~\ref{sec:ODE_absorb} contains the main results of this
paper. We show that $T_N$ is bounded by $Nt_N+O(N)$ and derive an
asymptotic development of $t_N$.  Section~\ref{sec:bounds} establishes a
explicit link between $T_N$ and the average completion time of one
algorithm. Finally, we show how this can be used to compute the completion
time of randomized algorithms in Section~\ref{sec:randomized_algo}.

\section{Formal description and notations}
\label{sec:notations}

Let $Y$ be a Markov chain on a finite state space $\Sp=\{0\dots S\}$ and
let denote $P$ its transition matrix. We assume that this chain has an
absorbing state, denoted $0$, all the other states being
transient.
If $I$ denotes the identity matrix, them $P-I$ can be written as:
\begin{equation*}
  P-I =\left[\begin{array}{cc}
      0 & \mathbf{0}\\
      Q^0 & Q
    \end{array}
  \right]
\end{equation*}
where $Q$ is a non-singular matrix such that for all $i,j$: $Q_{ii}<0$,
$Q_{ij}\ge0$ and $\sum_j Q_{ij}\le0$.  $Q^0$ is a vector such that for all
$i$, $\sum_j Q_{ij}+Q^0_{i}=0$.

We consider a Markov chain on $\Sp^N$ composed by $N$ copies of the
original Markov chain. Its state at time $t$ is denoted $\p{X_1(t),\dots
  X_N(t)}$. 
The evolution of the Markov chain is as follows:
\begin{enumerate}
\item at each time step, a chain $i\in\{1\dots N\}$ is picked uniformly at random
\item the $i$th chain $X_i(t)$ changes its state according to the
  transition matrix $P$. The states of the other chains do not change.
\end{enumerate}
Our goal is to compute the time for all chains to hit $0$ starting from a
state $(x_1\dots x_N)\in\Sp^N$. We define $T_N$ this hitting time:
\begin{equation}
  T_N \bydef \inf\left\{ t : \p{X_1(t),\dots X_N(t)}=(0\dots0) \right\}.
  \label{eq:T_N}
\end{equation}

\subsection{Notations}

For a state $x\in\Sp$, we denote by $\mathbf{e}_x$ the line vector that has
all coordinates equal to $0$ except for the $x$th one which is equal to
$1$. The vector $\mathbf{1}$ denotes the column vector with all coordinates
equal to $1$. For a line vector $\alpha$ and a matrix $P$, $\alpha P$
denotes the classical matrix product. For example, if $P$ is a $S\times S$
matrix, then $\alpha P \mathbf{1} = \sum_{i=1}^S\sum_{j=1}^S
\alpha_iP_{ij}$.

For each state $x\in\Sp$, we denote by $W(x)$ the hitting time of $0$
starting from $x$: if $Y$ is a Markov chain of transition probability $P$,
then:
\begin{equation*}
  W(x) \bydef \expect{ \inf\{ t : Y(t)=0 \} \mid Y(0)=x}.
\end{equation*}
Using the vector notation above, we have $\proba{\inf\{ t : Y(t)=0\} \ge
  1+i\mid Y(0)=x} = \mathbf{e}_xQ^{i}\mathbf{1}$ for all $i\in\N$. Thus,
$W$ can also be written
$ W(x) = \sum_{i=0}^\infty\mathbf{e}_xQ^i\mathbf{1}$.

\section{Hitting time and fluid approximation}
\label{sec:ODE_absorb}

In this part, we bound the expectation of $T_N$ using a deterministic
ordinary differential equation (ODE) approximation. In particular, we show
that $\expect{T_N}$ is bounded by $N$ times the time $t_N$ for the linear
ODE~\eqref{eq:ODE} to be at distance $1/N$ from its fixed point plus a
term of order $O(N)$. Moreover, this time $t_N$ is of order
$\Omega(\log(N))$, showing that the term in $O(N)$ becomes negligible
compared to $Nt_N$ as $N$ grows.

\subsection{A differential equation approximation}

For any state $x\in\Sp$ and any time step $k$, we define the quantity
$\bar{M}^N_x(k)$ to be the proportion of Markov chains that are in state
$x$ at time step $k$:
\begin{equation*}
  \bar{M}^N_x(k) = \frac1N\sum_{i=1}^N\mathbf{1}_{X_i(k) = x},
\end{equation*}
where $\mathbf{1}_{X_i(k)=x}$ equals $1$ if $X_i(k)=x$ and $0$ otherwise.
$\bar M^N(k)$ denotes the vector of all $\bar M^N_x(k)$ for $x\in\Sp$:
$\bar M^N(k)=\sum_{x\in\Sp} \bar M^N_x(k)\mathbf{e}_x$ where $\mathbf{e}_x$
denotes a unit vector having its $x$th coordinate equal to $1$ and the
others $0$. The process $\bar M^N(k)$ is a Markov chain: with probability
$\bar M^N_i(k)$, a chain that is in state $i$ is chosen and goes with
probability $P_{ij}$ in state $k$. This shows that the expected variation
of $\bar M^N(k)$ during one time step is:
\begin{align}
  \expect{\bar M^N(k+1) - \bar M^N(k) \mid \bar M^N(k)} &=
  \sum_{x\in\Sp}\sum_{j\ne i} \bar
  M^N_i(k)P_{ij}\frac{1}{N}(\mathbf{e}_j-\mathbf{e}_i)\nonumber\\
  &=\frac{1}{N}\bar M^N(k)Q.
  \label{eq:drift}
\end{align}
The function $f:m\mapsto mQ$ is called the \emph{drift} of the system.

Let us consider the system of differential equation corresponding to the
drift:
\begin{equation}
      \dot{m}(t) = m(t)\cdot Q.
  \label{eq:ODE}
\end{equation}

Equation~\eqref{eq:drift} shows that $M^N(k)$ can be described by a
stochastic approximation with constant step size $1/N$: it corresponds to a
Euler discretization of the ODE \eqref{eq:ODE} with a random noise $U$
(\emph{i.e.} such that $\expect{U^N(k+1) \mid \bar M^N(k)}=0$)
\begin{equation}
  \bar M^N(k+1) = \bar M^N(k) + \frac1N \lp \bar f(M^N(k)) + U^N(k+1) \rp.
  \label{eq:SA}
\end{equation}

Let us call $M^N(t)$ the state of the system when the time has been
rescaled by $t$: $M^N(t) = \bar M ^N(\lfloor tN \rfloor).$ Using classical
tools of stochastic approximation (Theorem~1 of \cite{benaimLeboudec} for
example), one can show that if $M^N(0)$ converges in probability to $m(0)$,
then $M^N(t)$ converges in probability to $m(t)$ uniformly on $[0;T]$:
\begin{equation*}
  \lim_{N\to\infty} \sup_{0\le t\le T} \norm{M^N(t)-m(t)} = 0
  \quad\mathrm{in~probability.} 
\end{equation*}

However, when one wants to compute the hitting time of the $0$ by $M^N(t)$,
this approximation is not accurate enough and leads to overestimated
bounds.  In the following, we will see how to link the hitting time of
$M^N(t)$ and the time for $m_0(t)$ to be greater than $1-1/N$.

\subsection{Link between the hitting time and the time to reach $1/N$}

Let us now look at the quantity $M^N_0(t)$ which is the proportion of
Markov chain in state $0$. The quantity $T_N$ can be defined as
\begin{equation*}
  T_N=\inf\{t:M^N_0(t)=1\}.
\end{equation*}
If $m(t)$ is the solution of the ODE~\eqref{eq:ODE}, one clearly has
$\lim_{N\to\infty}m_0(t)=1$ but unless if $m(t)$ starts exactly with
$m_0(t)=1$, the time to reach $m_0(t)=1$ is infinite:
$\inf\{t:m_0(t)=1\}=+\infty$.

Due to the discrete nature of $M^N_0(t)$, $M^N_0(t)$ takes values in
$\{0,\frac1N,\frac2N,\dots, \frac NN\}$. Thus, when $M^N_0(t)$ is greater
than $1-1/N$, it is equal to $1$. This suggests to introduce $t_N$, the
time for the ODE to be such that $m_0(t)\ge 1-1/N$:
\begin{equation}
  t_N\bydef\inf\{t: m_0(t)\ge 1-\frac1N\}.
  \label{eq:t_N}
\end{equation}

\begin{figure}[ht]
  \centering
  \subfigure[$N=20$]{\includegraphics[angle=-90,width=.48\linewidth]{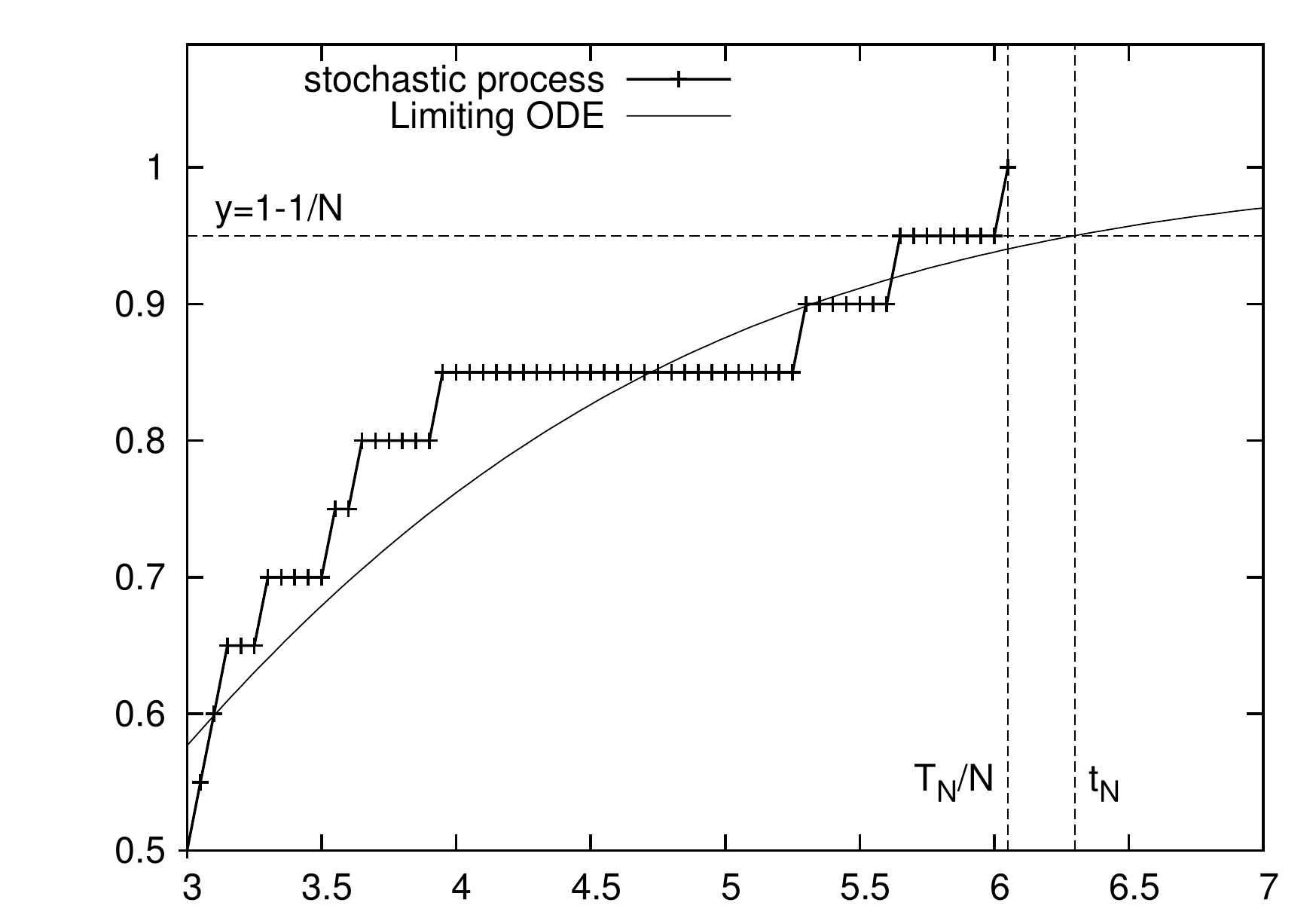}}
  \subfigure[$N=1000$]{\includegraphics[angle=-90,width=.48\linewidth]{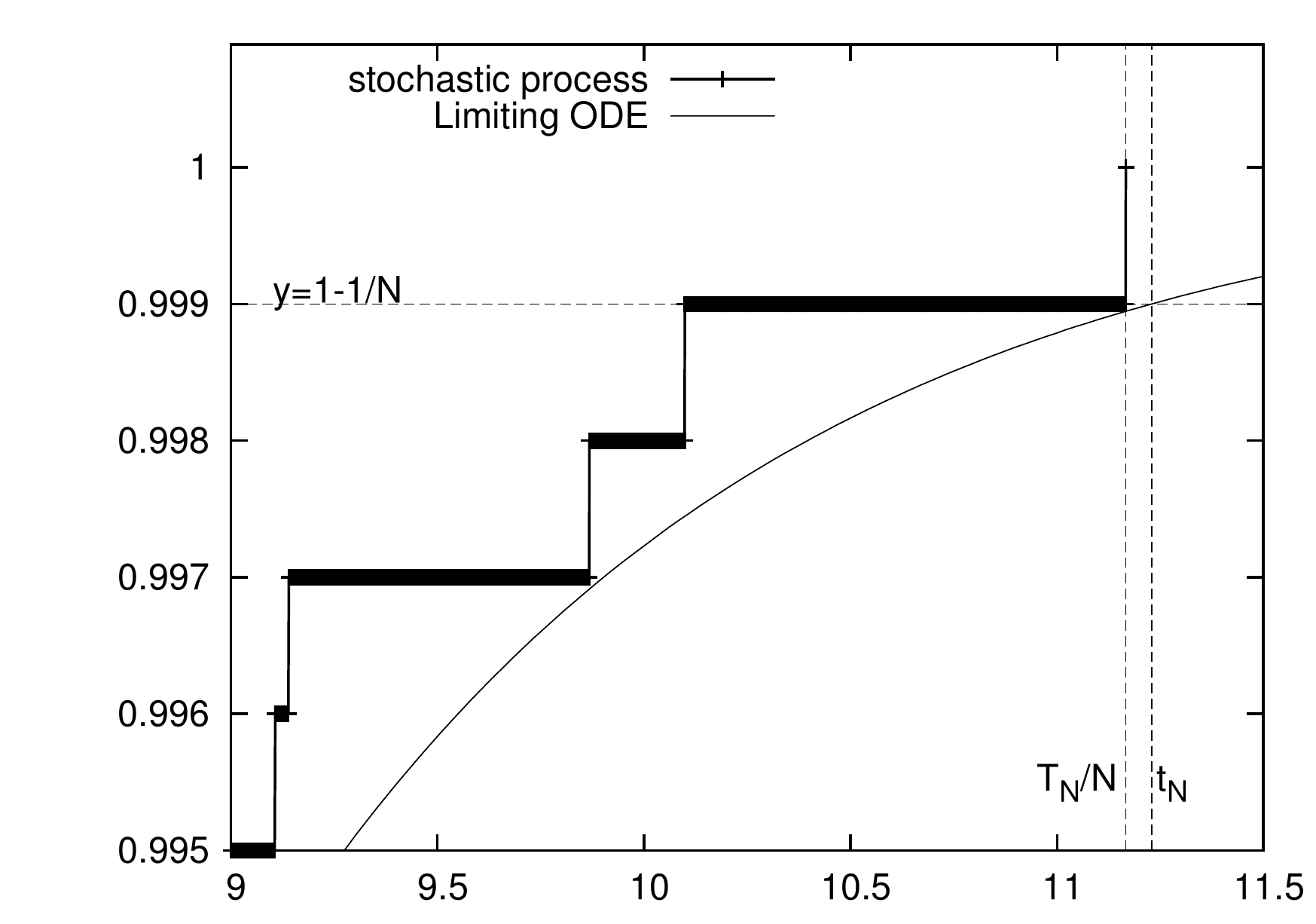}}
  \caption{Comparison of the hitting time of the stochastic of the system
    rescaled by $1/N$ and the time for the differential equation to reach
    $1-1/N$ for the coupon collector problem with $2$ cards. The smooth
    curve represents the differential equation $m_0(t)$, the dotted line is
    the line $1-1/N$ and the curve with the jumps represents $M^N(Nt)$ for
    one sample of the simulation for $N=20$ or $N=1000$.  For each curve,
    the hitting time of $1$ for the stochastic system is close to the
    hitting time of $1-1/N$ for the deterministic ODE. }
  \label{fig:simu_N=20}
\end{figure}

On Figure~\ref{fig:simu_N=20} are reported two simulations for the coupon
collector with $2$ cards. We compare the hitting time $t_N$ of the
stochastic system for $N=20$ and $N=1000$ with the time $t_N$ for the ODE
to reach $1-1/N$. The time of the stochastic system has been accelerated by
$N$. This suggests that $t_N$ is indeed a good estimate of $T_N/N$.

Classical stochastic approximation results show that the rate of
convergence of $M^N(t)$ to $m(t)$ is of oder $O(1/\sqrt{N})$. The bound is
too loose to guarantee the convergence of $t_N$ to $T_N$.  In the next
Theorem~\ref{th:ODE_absorb}, we use a slightly different approach to show
that $Nt_N$ is indeed a very good approximation of $T_N$.

\begin{theorem}
  \label{th:ODE_absorb}
  Let $t_N$ be defined by Equation~\eqref{eq:t_N} with $m$ satisfying the
  differential equation~\eqref{eq:ODE} with initial condition
  $m(0)=\alpha\bydef N^{-1}\sum_{i=1}^N\mathbf{e}_{x_i}$.  Then, hitting
  time $T_N$ of $(0\dots0)$ for the stochastic system composed of the $N$
  chains starting from $(x_1\dots x_n)$ satisfies:
  \begin{align*}
    \expect{T_N} &\le 
    N\p{t_N+\alpha(I-R)^{-1}\mathbf{1} +2\max_{j,k}Q^{-1}_{jk}} = Nt_N +
    O(N),
  \end{align*}
  where $R$ is the matrix defined by $R_{ii}=0$ and
  $R_{ij}=-\frac{Q_{ij}}{Q_{ii}}$ for $i\ne j$.
\end{theorem}

\begin{proof} The outline of the proof is as follows. The main idea is to
  write $T_N$ as the maximum of $N$ dependent random variables that
  correspond to the time for each chain to reach its absorbing state. Then
  we establish a relation between the expectation of this maximum and the
  tail behavior of the marginal distribution of each random variable. The
  marginal distribution for each chain follows a phase-type
  distribution. We show in Lemma~\ref{lem:discrete_vs_continuous} that its
  tail behavior can be approximated by the one of a continuous phase-type
  distribution that leads to the term $t_N$.

  \medskip
  
  Let us pick a chain $i\in\{1\dots N\}$ at random. The distribution of the
  initial state of $i$ is the distribution $\alpha$. If $i$ was alone, the
  probability for this chain to be in the absorbing state $0$ at time step
  $k$ starting from state $i$ would be $(P^{k})_{i,0} =
  \alpha(I+Q)^k\mathbf{1}$.

  When considering the system composed by the $N$ Markov chains, the Markov
  chain $i$ makes a transition with probability $1/N$. Thus, the
  probability for this particular chain to be in its absorbing state $0$ at
  time $k$ is $\alpha(I+N^{-1}Q)^k\mathbf{1}$. Therefore, the time $T_i$ at
  which the Markov chain $i$ has hit its absorbing state satisfies:
  \begin{equation}
    \proba{T^i\ge k} = \alpha(1+\frac{1}{N}Q)^k\mathbf{1}.
    \label{eq:distribution_i_k}
  \end{equation}

  If $i_1\dots i_N$ denotes a random permutation of $\{1\dots N\}$, then
  the time for all the Markov chains to have hit $0$ is $T_N=\max_{1\le
    k\le N}T^{i_k}$. The variables $T^{i_k}$ are identically distributed
  following the law given by Equation~\eqref{eq:distribution_i_k}. However,
  these variables are not independent.

  Using the union bound and the fact that $\proba{T_N\ge k}\le1$, we have: 
  \begin{equation*}
    \proba{T_N\ge k} \le \min\p{1,\sum_{k=1}^N\proba{T^{i_k}\ge k}} =
    \min\p{1,N\alpha(1+\frac{1}{N}Q)^k\mathbf{1}}. 
  \end{equation*}
  Therefore, the expectation of $T_N$ can be bounded by:
  \begin{align}
    \expect{T_N} = \sum_{k=1}^\infty\proba{T_N\ge k}
    &\le\sum_{k=1}^\infty\min\p{1,N\alpha(1+\frac{1}{N}Q)^k\mathbf{1}}\nonumber\\
    &\le \sum_{k=1}^{x_N-1}1 + \sum_{k=x_N}^\infty
    N\alpha(1+\frac{1}{N}Q)^k\mathbf{1} \label{eq:bound_x_N+}
  \end{align}
  where $x_N=\min\{k\in\N:\alpha(1+N^{-1}Q)^k\mathbf{1}\le
  2/N\}$. Moreover, using that
  $\sum_{k=x_N}^\infty(1+N^{-1}Q)^k=-(1+N^{-1}Q)^{x_N}NQ^{-1}$, we have
  \begin{align}
    \sum_{k=x_N}^\infty \alpha(1+\frac{1}{N}Q)^k\mathbf{1}
    &=-\alpha(1+\frac{1}{N}Q)^{x_N}NQ^{-1}\mathbf{1}\nonumber\\
    &\le \max_{j,k} \p{-Q^{-1}_{j,k}} N
    \sum_{j=1}^S\p{\alpha(1+\frac{1}{N}Q)^{x_N}}_{j}\nonumber\\
    &\le  \max_{j,k} \p{-Q^{-1}_{j,k}} N \cdot\frac2N\label{eq:comes_def_xn}\\
    &=2\max_{j,k} \p{-Q^{-1}_{j,k}},\nonumber
  \end{align}
  where the Inequality~\eqref{eq:comes_def_xn} comes from the definitions
  of $x_N$.

  Combining this inequality and the Equation~\eqref{eq:bound_x_N+}, we get:
  \begin{align*}
    \expect{T_N} \le x_N + 2N \max_{j,k} \p{-Q^{-1}_{j,k}}.
  \end{align*}
  The quantity $x_N$ is defined by
  $x_N=\min\{k:\alpha(1+\frac{1}{N}Q)^k\mathbf{1}\le2/N\}$. We show in
  Lemma~\ref{lem:discrete_vs_continuous} that:
  \begin{align*}
    x_N&\le N \inf\{t:\exp\p{tQ}\le\frac{1}{N}\}+
    \alpha(I-R)^{-1}\mathbf{1}\\
    &=N t_N + \alpha(I-R)^{-1}\mathbf{1},
  \end{align*}
  where $R$ is a matrix defined by $R_{ii}=0$ and
  $R_{ij}=-\frac{Q_{ij}}{Q_{ii}}$ for $i\ne j$.
\end{proof}

\subsection{Discrete and continuous phase-type distribution}

A random variable such that $\proba{X\ge t}=\alpha\exp(tQ)\mathbf{1}$ is
said to have a continuous phase-type distribution of parameter
$(Q,\alpha)$. Let $Y(t)$ be a Markov chain on $\Sp$ such that the rate of
transition from $i\ne0$ to $j\ne i$ is $Q_{ij}$ and the rate of transition
from $0$ to $i\ne0$ is zero. If $\alpha$ is the initial distribution of
$Y(0)$, then the time for $Y(.)$ to reach zero follows a phase-type
distribution of parameters $(Q,\alpha)$.  Similarly, a random variable such
that $\proba{X\ge k}=\alpha(1+Q/N)^k\mathbf{1}$ is said to have a discrete
phase-type distribution of parameter $(1+Q/N,\alpha)$. This corresponds to
the time for a discrete Markov chain of transition matrix $1+Q/N$ to reach
zero. We refer to \cite{latouche1999matrix}, Chapter~2 for more definitions
and properties of phase-type distributions.

The next lemma shows the relation between the tail of a continuous
phase-type distribution of parameter $(Q,\alpha)$ and the tail of a
discrete phase-type distribution of parameter $(1+Q/N,\alpha)$.

\begin{lemma}
  \label{lem:discrete_vs_continuous}
  Let $x_N=\min\{k\in\N:\alpha(1+N^{-1}Q)^k\mathbf{1}\le2/N\}$ and
  $t_N=N\min\{t\in\R:\alpha\exp(Qt)\mathbf{1}\le1/N\}$ be
  defined as in the proof of Theorem~\ref{th:ODE_absorb}.  Then:
  \begin{equation*}
    x_N \le N\p{t_N+\alpha(1-R)^{-1}\mathbf{1}},
  \end{equation*}
  where $R$ is a matrix defined by $R_{ii}=0$ and
  $R_{ij}=-Q_{ij}/Q_{ii}$ for $i\ne j$. 
\end{lemma}

\begin{proof}
  Let us consider the Markov chain $Y()$ associated with the continuous phase-type
  distribution of parameter $(Q,\alpha)$.
  With probability $\alpha_i$, the Markov chain starts in state $i$. If
  after $k$ jumps, the Markov chain is in state $i$, it stays there for a
  time $T_{ki}$ exponentially distributed of parameter $-Q_{ii}$ and then
  jump to a state $j\ne i$ with probability $-Q_{ij}/Q_{ii}$. Thus, the
  probability of being in state $i$ after $k$ jumps is $(\alpha R^k)_{i}$
  where $R$ denotes the matrix with $R_{ii}=0$ and $R_{ij}=-Q_{ij}/Q_{ii}$
  for $i\ne j$. Therefore, if $X^N$ is a continuous phase-type random
  variable of parameter $(\alpha,Q)$, $X^N$ has the same distribution as:
  \begin{equation}
    X = \sum_{k=0}^\infty \sum_{i=1}^S U_{ki} T_{ki},
    \label{eq:continuous_PT}
  \end{equation}
  where $U_{ki}$ are (dependent) Bernoulli variables of parameter $(\alpha
  R^k)_{i}$ and $T_{ki}$ are (independent) exponentially distributed
  variable of parameter $-Q_{ii}$.

  Similarly, the quantity $\alpha(1+N^{-1}Q)^k$ corresponds to the
  probability for a discrete phase-type random variable to be greater than
  $k$ and a variable $X$ with discrete phase-type distribution of parameter
  $(\alpha,I+Q/N)$ has the same distribution as:
  \begin{equation}
    \sum_{k=0}^\infty\sum_{i=1}^S U_{ki}T^{(N)}_{ki}
    \label{eq:discrete_PT}
  \end{equation}
  where $U_{ki}$ are the same variable as before and $T^{(N)}_{ki}$ are
  independent geometric random variables of parameter
  $-Q_{ii}/N$. 
  
  Since $T^{(N)}_{ki}$ is a geometric random variable of parameter
  $-Q_{ii}/N$, for all $t\in\R^+$, we have:
  \begin{align*}
    \proba{T^{(N)}_{ik}\ge tN} &= (1-\frac{Q_{ii}}{N})^{\ceil{tN}}
    \\ &
    \le (1-\frac{Q_{ii}}{N})^{tN+1}
    \\ &
    \le \exp(-Q_{ii}(t+\frac1N))
    \\&
    = \proba{T_{ki}\ge t+\frac1N}. 
  \end{align*}
  where the last inequality comes from the fact that $\log(1+x)\le x$. This
  shows that $N^{-1}T^{(N)}_{ki}$ is less than $T_{ki}+\frac1N$ (for the
  stochastic order). As pointed out, all the $T^{(N)}_{ki}$ and $T_{ki}$
  are independent in Equations \eqref{eq:continuous_PT} and
  \eqref{eq:discrete_PT}. Therefore, we can assume that
  $N^{-1}T^{(N)}_{ki}\le T_{ki}+\frac1N$ almost surely.  Using that for any
  positive random variable $A,B$ and any $t\in(0;\infty)$ and
  $\ell\in[0;t]$, we have:
  \begin{equation*}
    \proba{A+B\ge t} \le \proba{\p{A\ge\ell}\cup \p{B\ge t-\ell}}
    \le \proba{A\ge t-\ell} + \proba{B\ge \ell},
  \end{equation*}
  this shows that for any $\ell$:
  \begin{align}
    \proba{ \sum_{k=0}^\infty
      \sum_{i=1}^S U_{ki} T^{(N)}_{ki}\ge Nt}
    &\le\proba{ \sum_{k=0}^\infty \sum_{i=1}^S U_{ki} T_{ki}+\sum_{k=0}^\infty
      \sum_{i=1}^S U_{ki}\frac1N\ge t} \nonumber\\ 
    &\le\proba{\sum_{k=0}^\infty \sum_{i=1}^S U_{ki} T_{ki}\ge t-\ell}
    +\proba{\sum_{k=0}^\infty\sum_{i=1}^S U_{ki}\ge \ell}.
    \label{eq:sum_proba}
  \end{align}
  
  By Markov inequality, $\proba{\sum_{k=0}^\infty\sum_{i=1}^S U_{ki}\ge
    N\ell}\le\expect{\sum_{k=0}^\infty\sum_{i=1}^S U_{ki}}/(N\ell)$, with
  $\expect{\sum_{k=0}^\infty\sum_{i=1}^S U_{ki}}=\sum_k \alpha R^k
  \mathbf{1}=\alpha(1-R)^{-1}\mathbf{1}$. 
  This shows that if $\ell=\alpha(1-R)^{-1}\mathbf{1}$, then the second
  part of \eqref{eq:sum_proba} is less than $1/N$. Moreover, if
  $t=\ell+t_N$, the first part of \eqref{eq:sum_proba} is less than
  $1/N$. This shows that if $t\ge t_N + \alpha(1-R)^{-1}\mathbf{1}$,
  then:
  \begin{equation*}
    \proba{ \sum_{k=0}^\infty\sum_{i=1}^S U_{ki} T^{(N)}_{ki}\ge Nt}\le
    \frac2N.
  \end{equation*}
  Thus, this shows that $x_N\le N(t_N+\alpha(1-R)^{-1}\mathbf{1})$.

\end{proof}

\subsection{The logarithmic trend}

The quantity $\alpha\exp(Qt)\mathbf{1}$ is equal to one minus the
cumulative distribution function $F()$ of a continuous phase-type random
variable of parameter $(\alpha,Q)$. According to Theorem~2.7.2 of
\cite{latouche1999matrix}, there exist $\gamma>0$ and $k\ge0$ such that the
density of this variable $f()$ satisfies:
\begin{equation}
  f(t)=\gamma t^k \exp(-\nu t)+o(t^k \exp(-\nu t)),
  \label{eq:density_f}
\end{equation}
where $-\nu$ is a eigenvalue of $Q$ such that $\nu>0$ and $k+1\ge 1$ is the
multiplicity of the eigenvalue $-\nu$.

Equation~\eqref{eq:density_f} leads to the logarithmic bound for $T_N$,
expressed by the following theorem. 
\begin{theorem}
  \label{th:bound_log}
  Let $-\nu$ be the eigenvalue of $Q$ with the greatest real part and let
  $k+1\ge1$ denotes its multiplicity. Then, $\nu$ is real and positive and
  $T_N$ satisfies:
  \begin{equation*}
    \expect{T_N}\le \frac{1}{\nu} N\log(N) + \frac{k}{\nu} N\log\log N +O(N).
  \end{equation*}
\end{theorem}

\begin{proof}
  By Equation~\eqref{eq:density_f}, the cumulative distribution function
  $F$ satisfies:
  \begin{equation*}
    \alpha\exp(Qt)\mathbf{1}=1-F(t) = \int_{t}^\infty f(s)ds = \gamma t^k
    \exp(-\nu t)+o(t^k \exp(-\nu t)). 
  \end{equation*}
  Let $s_N(x)\bydef \nu^{-1}(\log(\gamma N)+k\log\log \gamma N-k\log
  \nu+x)$. For all fixed $x$, we have $s_N(x)=t_N+\nu^{-1}x+o(1)$.  Using
  that $\exp(-\nu s_N(x))=\frac{1}{\gamma N}(\log\gamma
  N)^{-k}\nu^k\exp(-x)$, the quantity $s_N(x)^k\exp(-\nu s_N(x))$ is equal
  to
\begin{align}
  v^{-k}(\log(\gamma N)+k\log\log\gamma N-k\log\nu+x)^k \frac{1}{\gamma N}(\log\gamma
  N)^{-k}\nu^k\exp(-x) \nonumber\\
  = \frac{1}{\gamma N}\exp(-x)(1+\frac{k\log\log\gamma N-k\log\nu+x}{\log
    \gamma N})^{-k}. 
  \label{eq:o(1)}
\end{align}  
The last factor of \eqref{eq:o(1)} goes to $1$ as $N$ goes to
infinity. Therefore, if $x>0$ (or $x<0$), then \eqref{eq:o(1)} is strictly
less (or greater) than $1/\gamma N$ if $N$ is large enough. This shows that
for all $\epsilon>0$, if $N$ is large enough, we have
\begin{equation*}
  1-F(s_N(-\epsilon)) < 1/N < 1-F(s_N(\epsilon)).
\end{equation*}

This shows that the number $t_N$ such that $t\ge t_N$ implies
$\alpha\exp(Qt)\mathbf{1}\le \frac1N$, is equal to:
\begin{equation}
  t_N=\frac{1}{\nu}\p{\log(\gamma N)+k\log\log N - k\log(\nu) + o(1)}.  
  \label{eq:bound_t_N}
\end{equation}
Combining \eqref{eq:bound_t_N} and Theorem~\ref{th:ODE_absorb} concludes the
proof of the theorem. 
\end{proof}

\subsection{Application to the coupon collector problem}
\label{sec:coupon_collector}

Let us consider the classical coupon collector problem: there are $N$
different types of coupon.  At each time step, a coupon of type $i$ is
picked at random where $i$ is drawn uniformly at random. It has been shown
in \cite{double-dixie} that the time to collect $T$ coupon of each type is
bounded by $N(\log N + (T-1)\log\log N +O(1))$. In this section, we show
that our approach allows one to retrieve this result directly.

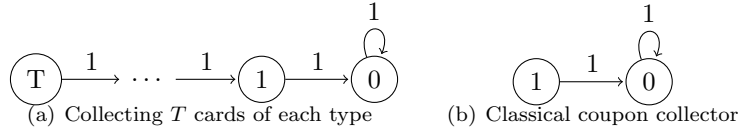
\begin{figure}[ht]
  \centering
  \begin{tabular}{ll}
    \subfigure[\label{fig:dixie}Collecting $T$ cards of each type]{
      \begin{tikzpicture}[shorten >=1pt,xscale=1.5]
        \node[circle,draw] at (0,0) (0) {0};
        \node[circle,draw] at (-1,0) (1) {1};
        \node at (-2,0) (2) {\dots};
        \node[circle,draw] at (-3,0) (T) {T};
        \draw (1) edge[->] node[above] {1} (0);
        \draw (2) edge[->] node[above] {1} (1);
        \draw (T) edge[->] node[above] {1} (2);
        \draw (0) edge[->,loop above] node[above] {1} (0);
      \end{tikzpicture}}&
    \subfigure[\label{fig:classical_coupon}Classical coupon collector]{
      ~~~~~~~~~\begin{tikzpicture}[shorten >=1pt,xscale=1.5]
        \node[circle,draw] at (0,0) (0) {0};
        \node[circle,draw] at (-1,0) (1) {1};
        \draw (1) edge[->] node[above] {1} (0);
        \draw (0) edge[->,loop above] node[above] {1} (0);
      \end{tikzpicture}~~~~~~~~~
    }
  \end{tabular}
  \caption{Markov chains used to represent the coupon collector problem:
    the state indicates the number of coupons of that remain to be
    collected. }
  \label{fig:coupon_dixie}
\end{figure}

Let us consider the Markov chain represented on Figure~\ref{fig:dixie}. Its
state space is $\{0\dots T\}$. The initial state is $T$ and for all $0<i\le
T$: $P_{i,i-1}=1$. The matrix $Q$ corresponding to this Markov chain is a
$T\times T$ matrix that has $-1$ on its diagonal and $1$ on its
sub-diagonal. On Figure~\ref{fig:classical_coupon} is represented the
particular case for $T=1$.

The ODE corresponding to this system is:
\begin{equation*}
  \left\{\begin{array}{lcll}
    \dot{m}_T(t) &=& -m_{T}(t) &\\
    \dot{m}_i(t) &=& -m_{i}(t)+m_{i+1}(t) & \quad \mathrm{for}~0<i<T
  \end{array}
  \right.
\end{equation*}
with $m_T(0)=1$ and $m_i(t)=0$ for $i\in\{0\dots T-1\}$.

A direct computation shows that $m_0(t)$ is the cumulative distribution
function of an Erlang variable of parameter $(T,1)$ (\emph{i.e.} the sum of
$T$ \emph{i.i.d.} exponential variable of parameter $1$) which can be
written: 
\begin{align*}
  m_0(t)&=1-\sum_{k=0}^{T-1}\exp(-x)x^k/k!\\
  &=1-\exp(-x)x^{T-1}/(T-1)!+O(\exp(-x)x^{T-2}).
\end{align*}
Using Theorem~\ref{th:ODE_absorb}, this shows that the time $T_N$ to
collect $T$ cards of each type is bounded by
\begin{equation*}
  \expect{T_N} \le \log N + (T-1)\log \log N + (T+2)N+O(1).
\end{equation*}
where $T+2$ comes from $Q^{-1}_{ij}=-\mathbf{1}_{j\le i}$ and $1-R=-Q$.

\section{Explicit formula for two particular cases}
\label{sec:bounds}

Theorem~\ref{th:bound_log} gives a precise idea on the behavior of $T_N$ in
the general case. However, the computation of the constants $\nu$, $k$ or
$O(N)$ can be difficult when the state space of the original Markov chain
is large. In this section, we derive explicit formulas for these constants
assuming that the hitting time for one Markov chain is bounded by $T$.

We first show that if the hitting time of the absorbing state is bounded by
$T$ for all single chain, then $T_N$ is less than $N^2T$
(Theorem~\ref{th:unbounded}), which is a loose bound in many cases.  When
the hitting time of the absorbing state is uniformly bounded by $T$ for all
initial states $x\in\Sp$, then $\expect{T_N}$ is bounded by $TN\log N+O(N)$
(Theorem~\ref{th:unif_bound}). 
At the end of the section, we provide two examples that shows that these
bounds are tight.

The results presented in this section remain valid if the state space of
the chain is countable instead of finite.

\subsection{Unbounded case}

If $W(x)$ denotes the expected hitting time of $0$ for a single Markov
chain, then the following results hold. 

\begin{theorem}
  \label{th:unbounded}
  The time $T_N$ such that all the chains have reached $0$ is bounded by:
  \begin{equation*}
    \expect{T_N}\le N\sum_{i=1}^NW(X_i(0))
  \end{equation*}
  In particular, if for all $i$, $W(X_i(0))\le T$, then  
  \begin{equation*}
    \expect{T_N} \le T N^2.
  \end{equation*}
\end{theorem}

\begin{proof}
  For all $i\in\{1\dots N\}$, let us call $R_i$ the time at which the
  Markov chain $i$ reaches $0$: $R_i=\inf\{t:X_i(t)=0\}$. It should be
  clear that $T_N=\max_{1\le i\le N}R_i\le \sum_{i=1}^NR_i$. Moreover, the
  hitting time for just one Markov chain when it performs one transition at
  each time is $W(X_i(0))$. As the probability for a Markov chain to
  perform one transition during one time step is $1/N$, we have
  $\expect{R_i}=NW(X_i(0))$.
\end{proof}

This result seems to be in contradiction with Theorem~\ref{th:bound_log}
that shows that if we fix a Markov chain, the expected hitting time of a
system composed of $N$ of these Markov chains is bounded by $O(N\log
N)$. However, the constant hidden in $O(N\log N)$ depends on the Markov
chain and the trend in $N\log N$ is only valid when the number $N$ goes to
infinity while the bound $TN^2$ only depends on $T$.  At the end of the
section, we provide an example that shows that this bound is tight. The
Markov chain used for this example depends on $N$.

\subsection{Uniformly bounded case $\sup_{x\in\Sp}W(x) \le T$} 

On Figure~\ref{fig:tight2}, we present a Markov chain that the bound of the
previous theorem is tight. This chain has a very particular shape: starting
from the initial state, the hitting time of the absorbing state is $1$ with
a probability $1-1/N^2$. With probability $1/N^2$, the chain jumps into a
state from which is takes $N^2(T-1)$ steps to hit $0$. This later causes
the hitting time to be large when multiple chain.

In this section, we show that if there are no such problematic states, the
bound on $T_N$ can be improved dramatically. More precisely, we show that
if the hitting time of $0$ is bounded by $T$ independently of the initial
state -- $\sup_{x\in\Sp} W(x)\le T$, then $T_N$ is of order $TN\log N$.
\begin{theorem}
  \label{th:unif_bound}
  If the hitting time for one chain is uniformly bounded by $T$
  (\emph{i.e.} $\sup_{x\in\Sp}W(x) \le T$), then the time $T_N$ such that
  all the chains have reached $0$ satisfies:
  \begin{equation*}
    \expect{T_N}\le TN \log N+ 2NT+1
  \end{equation*}
\end{theorem}

\begin{proof}
  Let $\F_t$ denotes the filtration associated to the process $X(t)$ and
  let us define the \emph{potential} of the system at time $t$, $\Phi_t$
  by:
  \begin{equation*}
    \Phi(t)\bydef\frac{1}{T}\sum_{i=1}^N W(X_i(t)).
  \end{equation*}
  $T_N$ is the time at which all $X_i(t)$ are equal to $0$ and can be
  written $T_N = \inf\{t:\Phi(t)=0\}$. 
  In the following, we first show that the time for $\Phi(t)$ to be lower
  than $1$ is less than $NT\log N + NT +1$ using
  \cite{tchibou-gast2010}. Then, we use Theorem~\ref{th:unbounded} to bound
  the remaining time by $NT$. 

  We say that a Markov chain is \emph{active} if it did not reach $0$. If
  an active Markov chain is picked in step $1$, then the potential will
  decrease in expectation by $1/T$. Let $\alpha(t)$ denotes the number of
  \emph{active} Markov chains at time $t$ (\emph{i.e.}
  $\alpha(t)=\sum_{i=1}^N \mathbf{1}_{X_i(t)\ne0}$). The probability of
  picking an active Markov chain is $\alpha(t)/N$. Therefore, the expected
  decrease of the potential between time $t$ and $t+1$ is:
  \begin{equation}
    \expect{\Phi(t+1) - \Phi(t)\mid \F_t} \le - \frac{\alpha(t)}{TN}.
    \label{eq:decrease_alpha}
  \end{equation}
  By hypothesis, $\sup_{x\in\Sp}W(x) \le T$. Thus, an active processor
  contributes at most $1$ to the potential and we have
  $\Phi(t)\le\alpha(t)$. Combining this with \eqref{eq:decrease_alpha}, we
  get:
  \begin{equation}
    \expect{\Phi(t+1)\mid \F_t} \le \Phi(t)\p{1-\frac{1}{NT}}.
    \label{eq:expected_decrease}
  \end{equation}
  
  Because of Equation~\eqref{eq:expected_decrease}, our potential function
  satisfies the hypothesis of Theorem~1 of \cite{tchibou-gast2010} with
  $m=1$ and $h(r)=1-1/(NT)$. According to this theorem, we have:
  \begin{align*}
    \inf\{t:\Phi(t)<1\} \le \lambda(\log(\Phi(0))+1)+1.
  \end{align*}
  where $\lambda=-1/\log(1-1/(TN))\le TN$ and $\Phi(0)\le N$. 

  By Theorem~\ref{th:unbounded}, when $T_N$ is less than $1$, $\sum_i
  W(X_i(t))\le T$. Therefore, the remaining time to hit $0$ is bounded by
  $NT$.
  
\end{proof}

\subsection{Comparison with previous bounds and tightness}
\label{sec:example_unbounded}

Theorems \ref{th:unbounded} and \ref{th:unif_bound} need stronger
assumptions than Theorems \ref{th:ODE_absorb} and \ref{th:bound_log} and
are often less precise.  However, their main advantage is to give explicit
formulas for the hitting time, even if computing the time $t_N$ or the
eigenvalue of the individual Markov chains is hard. This fact is important
in practical situation where the Markov chains often have a complicated
geometry. This is the case for the example of \cite{ludek} presented in the
next Section~\ref{sec:randomized_algo}.

The loss of precision of these bounds are well illustrated by the coupon
collector problem. Consider the Markov chain of Figure~\ref{fig:dixie} that
corresponds to the problem of collecting $T$ cards of each type. The
hitting time of $0$ from any state is clearly bounded by $T$. Therefore,
using Theorem~\ref{th:unif_bound}, one has $\expect{T_N}\le TN(\log
N+2)+1$. This bounds is worse by a factor $T$ compared with the bound
obtained by the ODE approach which was $N(\log N + (T-1)\log\log
N)+O(N)$. This is explained by the fact that Theorem~\ref{th:unif_bound}
does not take into account the particular shape of the Markov chain of
Figure~\ref{fig:dixie}: the bound of Equation~\eqref{eq:expected_decrease}
neglects the fact the hitting time starting from state $\{0\dots T-1\}$ is
strictly less than $T$.

\subsubsection{Tightness of the bounds of Theorems \ref{th:unbounded} and
  \ref{th:unif_bound}.}

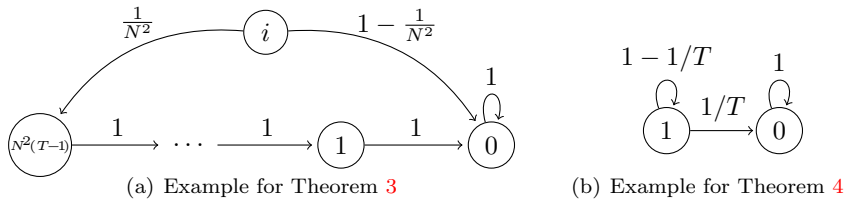
\begin{figure}[ht]
  \centering
  \begin{tabular}{ll}
    \subfigure[\label{fig:tight2} Example for  Theorem~\ref{th:unbounded}]{
      \begin{tikzpicture}[transform shape
        ,shorten >=2pt]
        \node[circle,draw] at (0,0) (0) {0};
        \node[circle,draw] at (-2,0) (1) {1};
        \node at (-4,0) (2) {\dots};
        \node[circle,draw,inner sep=0pt] at (-6,0) (T) {
          $\scriptscriptstyle N\!^2\!(T{\!-\!}1\!)$};
        \draw (1) edge[->] node[above] {1} (0);
        \draw (2) edge[->] node[above] {1} (1);
        \draw (T) edge[->] node[above] {1} (2);
        \draw (0) edge[->,loop above] node[above] {1} (0);
        \node[circle,draw] at (-3,1.5) (i) {$i$};
        \draw (i) edge[->,bend right] node[above] {$\frac{1}{N^2}$} (T);
        \draw (i) edge[->,bend left] node[above] {$1-\frac{1}{N^2}$} (0);
      \end{tikzpicture}
    }
    &
    \subfigure[\label{fig:tight1}Example for Theorem~\ref{th:unif_bound}]{
      ~~~~~~\begin{tikzpicture}[shorten >=1pt,xscale=1.5]
        \node[circle,draw] at (0,0) (0) {0};
        \node[circle,draw] at (-1,0) (1) {1};
        \draw (1) edge[->,loop above] node[above] {$1-1/T$} (1);
        \draw (1) edge[->] node[above] {$1/T$} (0);
        \draw (0) edge[->,loop above] node[above] {$1$} (0);
        \node at (0,-.5) {};
      \end{tikzpicture}~~~~~~
    } 
  \end{tabular}
  \caption{Two Markov chains used to show that the bounds of
    Theorems~\ref{th:unif_bound} and \ref{th:unbounded} are tight.}
  \label{fig:tight}
\end{figure}

A Markov chain that shows the tightness of the bound of
Theorem~\ref{th:unbounded} is represented on Figure~\ref{fig:tight2}. The
chain has $N+2$ states, denoted $\{0,\dots,N^2(T-1)\}\cup\{i\}$. Its
initial state is $i$. From $i$, the chain goes with probability $1/N^2$ to
state $N^2(T-1)$ and with probability $1/N^2$ in state $0$. From any state
$x\in\{1\dots N^2(T-1)\}$, the chain goes to state $x-1$ with probability
$1$.  For any state $x\in\{0\dots N^2 (T-1)\}$, the expected hitting time
of the state $0$ is $W(x)=x$. When starting in $i$, the expected hitting
time of $0$ is $T$. Therefore, Theorem~\ref{th:unbounded} shows that the
expected hitting time of $(0\dots0)$ for a system composed of $N$ of these
chains starting in $(i\dots i)$ is bounded by $TN^2$.

Let us compute a lower bound on the hitting time of $(0\dots0)$ starting
from $(i\dots i)$.  With probability $1-(1-1/N^2)^N$, there will be at
least one chain that needs $N^2(T-1)$ transitions to converges. At each
time step, this chain makes a transition with probability $1/N$. Thus, in
average this chain will take $(T-1)N^3$ time steps to converge. Since this
happens with probability $1-(1-1/N^2)^N$, a lower bound for the hitting
time of $(0\dots0)$ starting from everyone in state $i$ is
\begin{equation*}
  \expect{T_N}\ge N^3(T-1)(1-(1-1/N^2)^N) = (T-1)(N^2+o(1)).
\end{equation*}
This shows that the bound of Theorem~\ref{th:unbounded} is almost tight up
to an additive term of $N^2$.

To show that the bound of Theorem~\ref{th:unif_bound} cannot be improved
much without further assumption, let us consider the Markov chain
represented on Figure~\ref{fig:tight1}: it has two states and the
probability of going to $1$ to $0$ is $1/T$. The expected hitting time of
state $0$ starting from state $0$ is $T$. Theorem~\ref{th:unif_bound}
implies that $T_N\le NT(\log N+2)+1$. The exact value for $T_N$ is
$NT\sum_{i=1}^Ni^{-1}\approx NT\log N +\gamma NT + o(N)$ where
$\gamma\approx.57$ is the Euler–Mascheroni constant. This is close to our
theoretical bound up to an additive term of $(2-\gamma)NT$.

\section{Computing completion time of randomized algorithms}
\label{sec:randomized_algo}

In this section, we show how these results can be applied to study the
completion time of randomized algorithms and show how this can be used to
design efficient distributed protocols. 

\subsection{Completion time of randomized algorithm}

One motivation for this work comes from the study of the time for a set of
$N$ distributed randomized algorithm to all finish, in a scenario similar
to \cite{ludek}.

Let us consider that we want to solve a resource allocation problem among a
population of agent. We assume that we have a randomized algorithm that
converges to a stable allocation of the resource that is efficient but not
fair among different agents.  The final allocation might depend on the
random choices done by the algorithm and each allocation favors a different
group of agents. In order to improve the fairness of the equilibrium, we
consider the following scenario.  We execute $N$ independent copies of the
algorithm. At each time step, we do a step of computation of one algorithm
taken at random among the $N$ algorithms. After some time, the $N$
algorithms will have reached their stable allocation $S_1\dots S_N$. Since
each allocation is efficient, the resulting allocation will also be
efficient but the resulting allocation will be more fair since at each time
step, one allocation is picked at random among $S_1\dots S_N$.

If the original algorithm uses bounded memory, it can be represented by a
Markov chain with a finite state space. After some time, the Markov chain
will reach an absorbing state representing the fact that the algorithm
reached a stable allocation. The resulting algorithm can be represented by
$N$ independent Markov chains. At each time step, one Markov chain is
picked at random and performs one transition. Our framework, and in
particular Theorems~\ref{th:unif_bound} and \ref{th:unbounded} can be used
to compute the time to reach the final allocation if we know the time taken
by a single algorithm to converge.

\subsection{Correlated equilibria and distributed protocols}

These ideas are applied in \cite{ludek} to design a distributed algorithm
that converges to a fair and efficient allocation of wireless radio channel
to a set of user.

Their scenario is the following. There are $U$ users that want to share $C$
wireless channels.  The time is slotted and at each time slot, each user
can transmit data on one channel. It two or more users are transmitting on
the same channel at the same time, there are interferences and no data are
received. The only information available to a user before transmitting is
whether a given channel was used by one or more users at the previous time
slot.  Authors of \cite{ludek} proposed a distributed randomized algorithm
that converges to a constant assignment of the $C$ channels to $C$ of these
$U$ users while the others $U-C$ users do not transmit at all. This
algorithm guarantees a 100\% of utilization of the channels but is unfair
since $U-C$ users are not transmitting at all. The time of convergence of
the algorithm can be bounded by some constant $T$.

In order to improve the fairness of their algorithm, they introduced a
centralized entity that sends a correlation signals. At each time step $t$,
this entity sends a signal $n(t)\in\{1,\dots,N\}$ where $N$ is a predefined
constant (to be chosen by the entity). The signal $n(t)$ is picked
uniformly at random each time. The users keep $N$ copies of the previous
randomized algorithm. At time $t$, the users apply the algorithm number
$n(t)$. Even if this algorithm is no more completely distributed, it is
still scalable: the centralized entity has to broadcast a signal $n(t)$ at
time $t$ but it does not have to gather any information from the
agents. Moreover, this algorithm improves the fairness of the initial
algorithm.  As each algorithm is independent, if $N$ is large, each user
will be assigned in average to $C/U$ channels.  The more $N$ is large, the
more fair will be the allocation. However, a large $N$ slows the
convergence of the algorithm. An accurate bound on the time of convergence
of this algorithm allows one to choose the right compromise between the
speed of convergence and the performance of the algorithm. 

As the convergence time of each copy of the algorithm is bounded in
expectation by some $T$, this model satisfies the hypothesis of
Theorem~\ref{th:unif_bound}. Therefore, the convergence time of the whole
algorithm is bounded by $NT\log N + 2NT+1$.

\bibliographystyle{apt}
\bibliography{coupon_bibliography}

\end{document}